\documentclass[a4paper,12pt,oneside]{article} 
\title{\sc Theorems on Tangencies in Projective and Convex Geometry}
\author{{\sc Roland Abuaf}}
\date{\today}
\usepackage[dvips]{graphicx}

\usepackage{pstricks}

\DeclareGraphicsExtensions{eps}
\DeclareGraphicsExtensions{pdf}

\usepackage[T1]{fontenc}
\usepackage{amssymb}

\usepackage{amsmath}
\usepackage{latexsym}

\usepackage[latin1]{inputenc}

\newtheorem{theo}{Theorem}[subsection]

\newtheorem{prop}[theo]{Proposition}
\newtheorem{conjec}[theo]{Conjecture}
\newtheorem{defi}[theo]{Definition}
\newtheorem{nota}[theo]{Notations}

\newtheorem{conj}[theo]{Conjecture}

{
\noindent
\textit{\underline{Proof}} :\\
$\blacktriangleright\;$%
}
{\hspace{\stretch{1}}%
$\blacktriangleleft$}
{
\noindent
\textit{\underline{Preuve}} : (id\'ee)\\
$\blacktriangleright\;$%
}
{\hspace{\stretch{1}}%
$\blacktriangleleft$}

\newenvironment{skproof}
{
\noindent
\textit{\underline{Sketch of the proof}} :\\
$\blacktriangleright\;$%
}
{\hspace{\stretch{1}}%
$\blacktriangleleft$}

\begin{document}
\maketitle

\begin{abstract}

We discuss phenomena of tangency in Convex Optimization and Projective Geometry. Both theories have at disposal a powerful theory of duality. In both cases, the duality allows a nice interpretation of the contact locus of a hyperplane with an embedded variety. In this paper, we investigate more precisely some similarities between the theorems on tangencies existing in both theories. We focus in particular on a theorem of Anderson and Klee and its conjectural reformulation in Algebraic Geometry. If true, this conjecture would have significant consequences for Projective Geometry.

\end{abstract}

\newpage

\begin{section}{Introduction}

Let $X \subset \mathbb{R}^n$ be a compact convex body whose interior contains $0$. There have been considerable efforts to classify the \emph{singularities} of the points lying in the boundary of $X$. A clear picture of the situation was probably given for the first time by Anderson and Klee \cite{anderklee}.

\begin{defi} Let $X \subset \mathbb{R}^n$ be a compact convex body whose interior contains $0$, let $X^* \subset {\mathbb{R}^n}^*$ be its dual body and let $x \in \partial X$. We say that $x$ is an \emph{\textbf{$r$-singular point}} of $X$ if the \emph{\textbf{exposed face}} of $X^*$ relative to $x^{\bot}$ has dimension at least $r$.
\end{defi}

\begin{theo} 
Let $X \subset \mathbb{R}^n$ be a compact convex body whose interior contains $0$ and let $r \in \{0,...,n-1 \}$. The set of \emph{$r$-singular} points of $X$ can be covered by countably many compact subsets of finite $n-r-1$-dimensional Hausdorff measure.
\end{theo}
 
\bigskip

Now we turn to a similar situation in Projective Geometry. Let $X \subset \mathbb{P}_{\mathbb{C}}^n$ be an irreducible, non-degenerate projective variety. Zak found a bound for the dimension of the contact locus of any linear space with $X$ \cite{zak}.

\begin{theo} \label{zak} Let $X \subset \mathbb{P}_{\mathbb{C}}^n$ be an irreducible, non degenerate projective variety and let $L \subset \mathbb{P}_{\mathbb{C}}^N$ be a linear space. Denote by $X_L = \{x \in X,\, T_{X,x} \subset L \}$, we have the inequality:
$$ \dim(X_L) \leq \dim(L) - \dim(X)+b+1,$$
where $b = \dim X_{sing}$.
\end{theo}

Both theorems tell us that \emph{support loci} are subject to dimensional constraints. However, theorem 1.0.2 bounds the dimension of a family of hyperplanes when the dimension of the contact locus of the general member is known, whereas theorem 1.0.3 bounds the dimension of the contact locus of a single hyperplane. An Algebro-Geometric statement similar to theorem 1.0.2 was used in the first version of \cite{rasturm1}. Unfortunately, no convincing proof was given there.

\begin{conjec} Let $X \subset \mathbb{P}_{\mathbb{C}}^n$ be a non-degenerate, irreducible projective variety and let $X^* \subset \mathbb{P}_{\mathbb{C}}^N$ be its projective dual. Let $r \in \{ 0,...,n-1 \}$ and denote by $X^*\langle r \rangle = \{ H^{\bot} \in X^*, \, \dim \langle X_H \rangle \geq r \}$, where $\langle . \rangle$ denotes the scheme-theoretic linear span and $X_H$ is the tangency locus of $H$ with $X$. We have the inequality:
$$ \dim X^*\langle r \rangle \leq n-r-1.$$
\end{conjec}

\bigskip
 
I would like to thank Kristian Ranestad and Bernd Sturmfels for many interesting discussions on Convex Algebraic Geometry.

\end{section}

\begin{section}{Dualities and Contact Loci}
\begin{subsection}{A Common Setting for the Dualities}

Here we formulate, in a common language, the duality for convex bodies and for projective varieties. In the following, the space $\mathbb{E}^n$ either denotes the complex projective space $\mathbb{P}_{\mathbb{C}}^n$ or the real euclidean space $\mathbb{R}^n$. 
An object $X \subset \mathbb{E}^n$ refers to a compact convex body in $\mathbb{R}^n$ whose interior contains $0$ or to a reduced (irreducible) projective scheme in $\mathbb{P}_{\mathbb{C}}^n$. If $X$ is a convex body, then $\partial X$ is the boundary of $X$. If $X$ is a reduced projective scheme, then, by convention, $\partial X = X$. We denote by $\overline{Z}$ the convex hull or the Zariski closure of an object $Z \subset \mathbb{E}^n$.

\bigskip

\begin{defi}
Let $X \subset \mathbb{R}^n$ be a compact convex body whose interior contains $0$, let $y \in X$ and let $H \subset \mathbb{R}^n$ be a hyperplane. We say that $H$ has \emph{\textbf{contact}} with $X$ at $y$, if for all $x \in X$ we have $ \langle H^{\bot},x \rangle \leq 1$, and $\langle H^{\bot},y \rangle = 1$, where $\langle.,. \rangle$ is the evaluation pairing between ${\mathbb{R}^n}^*$ and $\mathbb{R}^n$.

\end{defi}
Note that if $H$ has contact with $X$ at $y$, then necessarily $y \in \partial X$.

\begin{defi}
Let $X \subset \mathbb{P}_{\mathbb{C}}^n$ be a reduced (irreducible) projective scheme and let $H \subset \mathbb{P}_{\mathbb{C}}^n$ be a hyperplane.

Let $y \in X_{smooth}$. We say that $H$ has \emph{\textbf{contact}} with $X$ at $y$ if $T_{X,y} \subset H$.

Let $y \in X_{sing}$. We say that $H$ has \emph{\textbf {contact}} with $X$ at $y$ if there exist sequences $(y_m) \in X_{smooth}$ and ${({H_m}^{\bot}) \in \mathbb{P}_{\mathbb{C}}^n}^*$ such that $H^{\bot} = \operatorname{lim} H_m^{\bot}$, $y = \operatorname{lim} y_m$ and $H_m$ has contact with $X$ at $y_m$ for all $m \in \mathbb{N}$.
\end{defi}

Now we can state both dualities in a common setting.
\newpage

\begin{theo}[Duality] Let $X \subset \mathbb{E}^n$ be an object. Consider the incidence $I_{X} = \{ (H^{\bot},x) \in {\mathbb{E}^n}^* \times \partial X, H \,\ has \,\ contact\,\ with\,\ X\,\ at \,\ x \},$ and the natural diagram:
\begin{figure}[!h]
\centering
\label{conormal}

\begin{tabular}{rrcll}
& & $I_X$ & & \\
& $\stackrel{q}\swarrow $ & & $\stackrel{p}\searrow$ & \\
$ {\mathbb{E}^n}^* \supset \partial (X^*) $ & & & & $ \partial X \subset \mathbb{E}^n$\\
\end{tabular}
\caption{conormal diagram}
\end{figure}

Let $X^* = \overline{q(I_X)}$. We have $I_{X^*} = I_{X}$. As a consequence, we have $X^{**} = X$ and $q(I_X) = \partial(X^*)$.

\end{theo}

Note that if $X \subset \mathbb{E}^n$ is a reduced projective scheme, then $q(I_X)$ is obviously Zariski closed. In this case, the equality $q(I_X) = \partial(X^*)$ is trivial since, in our notations, $\partial X^* = X^*$. Note also that, by construction, for all $H^{\bot} \in q(I_X)$, we have:
$$ x \in p(q^{-1}(H^{\bot})) \Leftrightarrow H \,\ has \,\ contact \,\ with \,\ X \,\ at \,\ x.$$ The object $p(q^{-1}(H^{\bot}))$ is called the \emph{\textbf{contact locus}} of $H$ along $X$. In Convex Geometry, the set $p(q^{-1}(H^{\bot}))$ is often called the \emph{exposed face} of $X$ relative to $H$, while in Projective Geometry it is known as the \emph{tangency locus} of $H$ along $X$. The duality says that the set of hyperplanes which have contact with $X$ at $x$ is equal to the contact locus of $x^{\bot}$ along $X^*$. That is, for all $x \in X$, we have:
$$ H^{\bot} \in q(p^{-1}(x)) \Leftrightarrow x^{\bot} \,\ has \,\ contact \,\ with \,\ X^* \,\ at \,\ H^{\bot}.$$

\bigskip

A proof of this result can be given in the framework of Lagrangian Geometry. Indeed the notion of \emph{contact} allows one to define in a uniform way Lagrangian manifolds for Projective and Convex Geometry. 
\end{subsection}

\begin{subsection}{The Principle of Anderson and Klee} 
In this section, we formulate the principle of Anderson and Klee in a common setting for Projective Geometry and Convex Geometry.

\begin{nota} Let $X \subset \mathbb{E}^n$ be an object. The \emph{\textbf{linear span}} of $X$, which we denote by $\langle X \rangle$ is the smallest linear subspace of $\mathbb{E}^n$ which contains $X$.
\end{nota}

In the case $Z \subset \mathbb{E}^n$ is a non-reduced scheme, the subspace $\langle Z \rangle $ is the scheme-theoretic linear span of $Z$.

\begin{defi} Let $X \subset \mathbb{E}^n$ be an object. A point $x \in X$ is said to be a \emph{\textbf{$r$-singular point}} in $X$ if $\dim \langle q(p^{-1}(x)) \rangle \geq r$. The set of $r$-singular points of $X$ is denoted by $X\langle r \rangle$.
\end{defi}

The following result is the archetype of the theorem on tangencies which should be true in all geometries. It was proven by Anderson and Klee (see \cite{anderklee}, or \cite{schnei} for a modern presentation) in the context of Convex Geometry and was used in the first version of \cite{rasturm1} in the context of Projective Geometry.

\begin{conj}
Let $X \subset \mathbb{E}^n$ be an object, we have the inequality:
$$ \dim X\langle r \rangle \leq n-r-1.$$
\end{conj}
Here the dimension must be understood as the Hausdorff dimension or the algebraic dimension, depending on the theory. Note that this conjecture is of cohomolgical nature in Algebraic Geometry. A similar statement can be formulated as follows.

\begin{conj}  \label{main} Let $X \subset \mathbb{P}^n$ be a smooth, irreducible, non-degenerate projective variety. Denote by $X^*\langle r \rangle$ the set:
$$ X^*\langle r \rangle = \{H^{\bot} \in X^*, \,\, \dim H^{0}(\mathcal{J}_{p(q^{-1}(H^{\bot}))}(1)) \leq n-r \}.$$ We have the inequality:
$$ \dim X^* \langle r \rangle \leq n-r-1,$$ for all $r \in \{0,...,n\}.$
\end{conj}

Here $\mathcal{J}_{p(q^{-1}(H^{\bot}))}$ denotes the ideal sheaf of $p(q^{-1}(H^{\bot})$ in $\mathbb{P}^n$, where $p,q$ are defined in the conormal diagram. One expects that similar bounds could be found for the dimension of the set of points $H^{\bot} \in X^*$ such that $\dim H^{0}(\mathcal{J}_{p(q^{-1}(H^{\bot}))}(k))$ is "small" enough.
 
\bigskip

Using the theory developped by Hironaka around the notion of normal flatness \cite{hiro} and a result of L\^e-Teissier \cite{lete}, one can prove the following result in Projective Geometry.

\begin{prop} \label{reduced} Let $X \subset \mathbb{P}^n$ be an irreducible, reduced, non-degenerate projective variety. Let 
$$\widetilde{X}^*\langle r \rangle = \{ H^{\bot} \in X^*, \,\, \dim H^{0}(\mathcal{J}_{|p(q^{-1}(H^{\bot}))|_{red}}(1)) \leq n-r \}.$$ We have the inequality:
$$ \dim \widetilde{X}^*\langle r \rangle \leq n-r-1.$$
\end{prop}
Here $|p(q^{-1}(H^{\bot}))|_{red}$ is the reduced space underlying $p(q^{-1}(H^{\bot}))$.
\end{subsection}
\end{section}

\begin{section}{Applications to Projective Geometry}
If true, conjecture \ref{main} would have significant consequences for Projective Geometry. In fact, even proposition \ref{reduced} can be used to prove a generalization of Severi's theorem.

\begin{nota} Let $X \subset \mathbb{P}_{\mathbb{C}}^n$ be an irreducible projective variety. We denote by $X^*(r)$ the set $X^*(r) = \{ H \in X^*, \,\, \dim p(q^{-1}(H)) \geq r \}.$
\end{nota} 

\begin{theo} \label{severi} Let $X \subset \mathbb{P}_{\mathbb{C}}^5$ be a smooth, irreducible, non-degenerate projective surface and let $X^* \subset {\mathbb{P}^5}^*$ its projective dual. We have $\dim X^*(1) \leq 2$, with equality if and only if $X$ is the Veronese surface.
\end{theo}

\begin{skproof} By assumption $X \neq \mathbb{P}^2$, so $X$ does not contain a 2-dimensional family of lines. As a consequence of proposition \ref{reduced}, we see that $\dim X^*(1) \leq 2$.

\smallskip
 Assume that $\dim X^*(1) = 2$, proposition \ref{reduced} again shows that for all $H^{\bot} \in X^*(1)$, the curve-components of $(H \cap X)_{sing}$ are plane curves.

Let $H^{\bot} \in X^*(1)$ be a general point and let $k$ be the maximum of the degree of the curve-components of $|(H \cap X)_{sing}|_{red}$. Assume that $k \geq 3$. Then, there is a plane curve, say $C$, in $|(H \cap X)_{sing}|_{red}$ such that all lines in $\langle C \rangle $ are trisecants to $X$. But this is true for general $H^{\bot} \in X^*(1)$, so that a careful count of dimension shows that we have a 4-dimensional family of trisecants to $X$. This is impossible by the trisecants lemma.

\smallskip

As a consequence, the smooth surface $X$ is covered by a $2$-dimensional family of conics, it is the Veronese surface.

\end{skproof}

Note that theorem \ref{severi} obviously implies Severi's original result. Indeed, if $X \subset \mathbb{P}_{\mathbb{C}}^5$ is a smooth, irreducible, non-degenerate surface whose secant variety does not cover the ambiant space, then Terracini's lemma implies that $\dim X^*(1) = 2$. Another proof of Severi's result, relying on similar techniques as the above ones, is due to Zak and is a consequence of theorem \ref{zak}. Hence, one may hope that theorem \ref{zak} and conjecture \ref{main} could be considered in a common setting. As such, these results are perhaps incarnations of a deeper principle, which has yet to be discovered.

\end{section}

\newpage

\bibliographystyle{alpha}

\bibliography{bibli}

\end{document}